\documentclass{article}[11pt]

\usepackage{amssymb,amsfonts,amsmath,mathtext}
\usepackage{cite,enumerate,float,indentfirst}
\usepackage{wrapfig}
\usepackage{graphicx}
\usepackage{subcaption}

\title{Rectangular knot diagrams classification with deep learning}

\author{L.H. Kauffman 
\thanks{Department of Mathematics, Statistics and Computer Science
851 South Morgan Street
University of Illinois at Chicago
Chicago, Illinois 60607-7045
and
Department of Mechanics and Mathematics,
Novosibirsk State University,
Novosibirsk
Russia}
\thanks{Kauffman's work was supported by the 
Laboratory of Topology and Dynamics, 
Novosibirsk State University 
(contract no. 14.Y26.31.0025 with the Ministry of Education and Science of the Russian Federation).}
\and
N.E. Russkikh
\thanks{AcademGene LLC, 630090 Novosibirsk, Russia, and  A.P. Ershov Institute of Informatics Systems, 630090 Novosibirsk, Russia}
\and 
I.A. Taimanov
\thanks{Novosibirsk State University, 630090 Novosibirsk, Russia, and Sobolev Institute of Mathematics, 630090 Novosibirsk, Russia;
e-mail: taimanov@math.nsc.ru.} 
\thanks{Taimanov's work was supported by the Laboratory of Topology and Dynamics, 
Novosibirsk State University (contract no. 14.Y26.31.0025 with the Ministry of Education and Science 
of the Russian Federation).}
\thanks{The preliminary results of our research were presented by N.E. Russkikh in the talk of ``Recognition of the unknot via neural networks (Preliminary results)'' at ``December Readings in Tomsk - 2019'' (December 13, 2019). We thank participants of this conference for helpful discussions.}}

\begin{document}

\date{}
\maketitle

\begin{abstract}
In this article we discuss applications of neural networks to recognising knots and, in particular, to the unknotting problem. One of motivations for this study is to understand how neural networks work on the example of a problem for which rigorous mathematical algorithms for its solution are known. We represent knots by rectangular Dynnikov diagrams and apply neural networks to distinguish a given diagram’s class from the given finite families of topological types. The data presented to the program is generated by applying  Dynnikov moves to initial samples. The significance of using these diagrams and moves is that in this context the problem of determining whether a diagram is unknotted is a finite search of a bounded combinatorial space. 
\end{abstract}

\section{Introduction}

In this article we discuss the applications of neural networks to recognising knots and, in particular, to the unknotting problem which appears to be the most interesting application.

We represent knots by rectangular diagrams (see the text of the paper for the definition of these diagrams and for the moves that we use in transforming the diagrams)  and apply networks to distinguish a given diagram’s class from the given finite family of classes (topological types of knots which include the unknot class) that  may be represented by the diagram.

The data we present to our program is self-generated in the sense that we begin with samples of knots and unknots and then apply the Dynnikov moves (see the text of the paper for a definition of these moves)   to produce the larger collections that the program uses. We have two types of data. In one case we use internal moves only. In the second case we use external moves as well. In the second case we find that the program has a significantly harder time identifying the data. This provides us with a material difference that can be used for the next stages of the research. The significance of using rectangular diagrams and the Dynnikov moves is that in this context the problem of determining whether a diagram is unknotted is a finite search of a bounded combinatorial space. This finiteness is not available for algorithms using the Reidemeister moves or other formulations of the knot theory.

In comparison to other problems to which neural networks are applied the knot recognition problem has two special features:

\begin{itemize}
    \item 
there is no a natural data set and, in particular, the efficiency of networks is evaluated on generated data sets and strongly depends on generation methods;

\item
there are several mathematical algorithms for solving the unknotting problem (see Section 2). The comparison of them with the applied neural networks may help in understanding how these neural networks function. 
\end{itemize}

\section{Algorithms in three-dimensional topology}

The most famous recognition problem of three-dimensional topology which was solved algorithmically is the unknotting problem for knots. It has a long story. 

\begin{itemize}
\item
The first solution to the unknotting problem was found by Haken \cite{Haken}
who introduced for that the theory of normal surfaces in three-manifolds. 
\end{itemize}

Although his algorithm was quite inefficient,  this approach led to a special study of Haken, or sufficiently large, three-manifolds and revolutionized three-dimensional topology.
The idea consisted in finding an embedded  two-dimensional disk which is bounded by the knot. By the Papakiriakopolous theorem, such a disk exists if and only the given knot is the unknot \cite{Pap} .  

\begin{itemize}
\item
The algorithm by Birman and Hirsch used so-called braid foliations and also consists in looking for a disk bounded by a knot \cite{BH}. Their considerations are based on 
normal surface theory; 

\item
A different approach was proposed by Dynnikov who proved that 
any presentation of the unknot by a rectangular diagram admits a monotonic simplification by elementary moves \cite{Dynnikov}. 
\end{itemize}

The work by Dynnikov was preceded by another: his approach based on Gaussian codes of knots which led to a partial algorithm that works very fast for some quite complicated representations. However it is not effectively applied to all representations of the unknot. This
partial algorithm was realized on computers \cite{ADP} and this software was used for experiments
with different representations.

In our research we also use rectangular diagrams.

Two other solutions of the unknotting problem are based on algebraic invariants of low-dimensional manifolds:

\begin{itemize}
\item
Knot Floer homology contains information on the Seifert genus of a knot, which is the minimal genus of embedded surfaces bounded by the knot. By the Papakyriakopoulos theorem, the Seifert genus vanishes exactly for the unknot. The combinatorial representation of knot Floer homology gives an algorithm for computing the Seifert genus, and hence for the recognition of the unknot \cite{MOS};

\item
Kronheimer and Mrowka proved that the property that the reduced Hovanov cohomology of a knot has rank one distinguishes the unknot \cite{KM}.
\end{itemize}

Speaking about the algorithms, we refer to the complexity of this recognition problem.
By using normal surface theory,   Hass, Lagarias, and Pippenger showed that unknottedness
is in the complexity class NP \cite{HLP}.
On the other hand, it was proved by Kuperberg that, assuming the generalized Riemann hypothesis, 
the knottedness is also in the NP class, which implies that unknottedness is in co-NP  \cite{Kuperberg}.

A fundamental problem of knot theory reads

\vskip2mm

{\sc Is the Jones polynomial of a knot $K$ equal to one, $J(K)=1$, if and only if $K$ is the unknot? }

\vskip2mm

The positive answer to this question would give another solution of the unknottedness problem.

For general knots and links the problem of algorithmic recognition was positively solved by Matveev \cite{Matveev}
who actually established this result for more general class of Haken manifolds.

One of the traditional approaches for comparing general links consists in comparing their planar diagrams. It is know that two planar diagrams define the same link if and only if they can be related by a sequence of so-called Reidemeister moves. 

Few years ago 
\begin{itemize}
    \item 
    Coward and Lackenby proved that there exists a superexponential function $f(n_1,n_2)$ such that for any two planar diagrams of a link with $n_1$ and $n_2$ crossings there are at most $f(n_1,n_2)$ Reidemeister moves which take one diagram into another \cite{CL};
    
    \item
    Lackenby proved that a planar diagram of the unknot with $n$ crossings can be converted into the diagram with no crossings by at most $(236n)^{11}$ Reidemeister moves \cite{L}.
\end{itemize}

We remark that the theorem by Coward and Lackenby gives a new solution of the algorithmic recognition problem for links.

For a recent survey of algorithmic problems of three-dimensional topology we refer to \cite{L2020}.

\section{Rectangular diagrams}
\subsection{Definition} \label{rect_definition}
\begin{figure}[H]
  \center
  \includegraphics[width=0.99\linewidth]{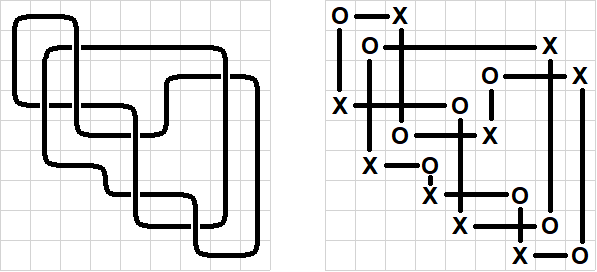}
  \caption{An example of the knot with corresponding rectangular diagram}
  \label{fig:init_diag}
\end{figure}

To apply computer science methods to mathematical objects, one must find computer-compatible representation for them. In our case, we use rectangular diagrams \cite{Dynnikov}. Basically, any knot diagram can be transformed in such way that all edges are either vertical or horizontal. Less obvious, we can further transform it so all vertical lines are overcrossing lines. Finally, we place all vertices on the grid in such way that only one edge belongs to each horizontal or vertical line of the grid. An example of the rectangular diagram  is depicted at figure [\ref{fig:init_diag}]. 

At this point, we're ready to represent the diagram numerically. Namely, we label each vertex with 
$X$ or $O$, consecutively changing them, and after that encode each horizontal edge with coordinates of $X$ and $O$ end. Since on each vertical line is only one edge, the diagram is represented by two permutations - coordinates of all $X$ and $O$ ends. There are two permutations, one corresponding to $X$ and the other corresponding to $O$. The $X$ permutation is a vector list of the columns that appear with an $X$ as one goes through the rows in order from top to bottom. The $O$ permutation is a vector list of the columns labeled $O$ as one goes through the rows in order from top to bottom. We indicate the two permutations by a matrix of two rows with the $X$ permutation preceding the $O$ permutation, which is a numerical  representation of the diagram. As an example, the diagram from the Fig.1, corresponds to the following representation:
$$
\begin{pmatrix}
3 & 8 & 9 & 1 & 6 & 2 & 4 & 5 & 7 \\
1 & 2 & 6 & 5 & 3 & 4 & 7 & 8 & 9
\end{pmatrix}, 	
$$
where the top row corresponds to $X$-vertices, and bottom - to the $O$ ones.
The \textbf{complexity} of the rectangular diagram is defined as the number of horizontal edges.

\subsection{Elementary moves}
\begin{figure}[H]
  \center
  \includegraphics[width=0.9\linewidth]{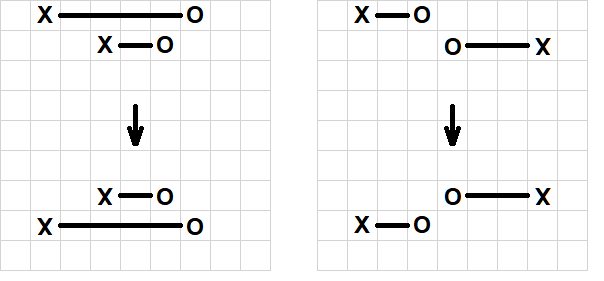}
  \caption{An example of the internal switch move.}
  \label{fig:move1_diag}
\end{figure}

\begin{figure}[H]
  \center
  \includegraphics[width=0.5\linewidth]{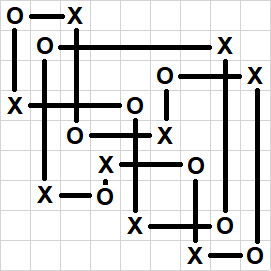}
  \caption{A diagram obtained from the initial one with the second Dynnikov move.}
  \label{fig:move2_diag}
\end{figure}

\begin{figure}[H]
  \center
  \includegraphics[width=0.5\linewidth]{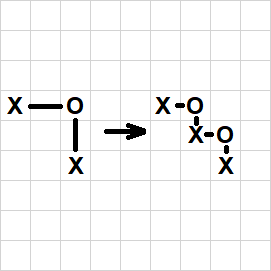}
  \caption{An example of the stabilization move}
  \label{fig:move_stab}
\end{figure}

There are three kinds of class-preserving rectangular diagram moves, which called Dynnikov moves:
\begin{itemize}
\item Neighboring edges can be switched if they are not interleaving. An example of that is depicted at Fig.[\ref{fig:move1_diag}]. A diagram which is one move away from the initial one is depicted at \ref{fig:inter_switch_diag}. Please note that switched edges share an endpoint, and despite that they're not interleaving.
\item The top and the bottom (the leftmost and the rightmost) edges can be switched under the same assumptions (it is called external switch, or translation move). The example can be seen at the \ref{fig:move2_diag}.
\item And finally, corners can be eliminated (or created) the way is depicted at the Fig. [\ref{fig:move_stab}]. This move is called (de)stabilization. 
\end{itemize}

Performing these moves in terms of permutation coding is described in the Appendix [\ref{section: perm}].

This set of moves is worth considering because it is proved in \cite{Dynnikov} that any unknot diagram can be monotonically (in terms of complexity) simplified to the diagram with no crossings. This means that determining whether a diagram is unknotted is a finite search using this notion of complexity.

\section{Application of deep learning to the unknottedness problem}
Despite that mathematical community was mostly focused on the algorithms which simplify given knot diagram, we take the different path.
Our program answers the following question(which extends unknottedness problem): to which knot class does (up to 36 classes with less than or equal to 8 crossings in the minimal diagram) a given rectangular diagram belongs? We  answer this question considering the diagram as-is, without attempts of its simplification. 
We use a collection of knot types from the tables of knots up to eight crossings. The data we present to our program is self-generated in the sense that we begin with samples of knots and unknots and then apply the Dynnikov moves to these to produce the larger collections that the program uses. We have two types of data. In one case we use internal moves only. In the second case we use external moves as well. In the second case we find that the program has a signifcantly harder time identifying the data. This provides us with a material difference that can be used for the next stages of the research.

 In terms of deep learning, our program performs multiclass classification task.  Deep learning approach was used to address this task in \cite{Vand} though they used a different knot parametrization: they represented a knot as a polymer conformation, consecutively attached unit length rods. Their sampling procedure is also different: they always start with an unknot (random circular conformation) and evolve it to the final configuration possibly changing its class. The class is derived afterwards using some knot invariants. Oppositely, we use label-preserving perturbations of the diagram, starting with the diagram with the minimal number of crossings. This lets us consider knots with arbitrary complexities and classes since we don't need any algorithms to assign labels for the training set.

\begin{tabular}{ |p{3cm}|p{3cm}|p{3cm}|  }
\hline
\multicolumn{3}{|c|}{Comparison of our approach with \cite{Vand}} \\
\hline
Property & Vandans Paper& Our approach\\
\hline
Number of classes & 5 & 36 \\
Architecture & BiLSTM   & BiLSTM \\
Labeling & Alexander polynomial & Label of the initial diagram \\

\hline
\end{tabular}

\section{Deep learning foundations}
\subsection{Network architecture} \label{architecture}
 We consider a knot as a sequence of coordinates of its horizontal edges in the rectangular diagram (see section \ref{rect_definition}). Discovering architectures of neural networks which fit best to sequential nature of the data is a long-standing field in the deep learning research. We chose the classic LSTM layer \cite{LSTM} to address it. We used a bidirectional variation of it  \cite{BiRNN}. The layer is composed of a \textbf{cell}  $c$ (the memory part of the LSTM unit) and three ``regulators'', usually called \textbf{gates}, of the flow of information inside the LSTM unit: an input gate $i$, an output gate $o$ and a forget gate $f$. All $c, h, i, f$ and $o$ are vectors of the same representation dimension $d$, which we have set to 1024.  Intuitively, the cell is responsible for keeping track of the dependencies between the elements in the input sequence. The input gate controls the extent to which a new value flows into the cell, the forget gate controls the extent to which a value remains in the cell and the output gate controls the extent to which the value in the cell is used to compute the output activation of the LSTM unit.  To sum up, the LSTM layer consumes an input sequence step by step (one step for each vertex $x_t \in \mathbf{R}^2$), updating  $h$ and $c$ via the following set of equations:
 
$$
f_t = \sigma(W_f x_t + U_f h_{t-1} +b_f),
$$

$$
i_t = \sigma(W_i x_t + U_i h_{t-1} +b_i),
$$

$$
o_t = \sigma(W_o x_t + U_o h_{t-1} +b_o),
$$

$$
c_t = f_t \odot c_{t-1} + i_t \odot tanh(W_c x_t + U_c h_{t-1} + b_c),
$$

$$
h_t = o_t \odot tanh(c_t),
$$
Matrices $W_f$, $W_i$, $W_o$, $W_c$ are called \textbf{weights}, vectors $b_f$, $b_i$, $b_o$ and $b_c$ are called biases, and their elements are learnable parameters which are adjusted during the training of the neural network.
Notation $\odot$ stands for elementwise multiplication and $\sigma$ denotes sigmoid activation function 
$$\sigma (z) = \dfrac{1}{1+\exp{(-z)}}$$ ,
which is also applied elementwise and squeezes each coordinate of a vector  into $(0,1)$ .

In order to obtain fixed size representation out of LSTM layer outputs, we perform an operation which is called \textbf{global average pooling}:
$$
\Tilde{h} = \dfrac{1}{N}\sum_{t=1}^{N} h_t
$$

Then, this fixed size representation is passed to the \textbf{fully connected layer}, which is, mathematically, is just a matrix multiplication:

$$
z = W_{fc}h,
$$
where $W_{fc}$ has the size of (36, 2048) since there are 36 classes of knots with crossing number up to 8 and we set the dimension of our LSTM hidden representation to 1024, and using bidirectional modification doubles the dimension by using representations for another direction. Finally, we apply a transformation to $z$ which is called a \textbf{softmax function}:

$$
\widetilde{y}_i = \dfrac{e^{z_i}}{\sum_{i=1}^{36}e^{z_i}}.
$$
Each $y_i$ corresponds to the model's confidence that the sample belongs to the class $i$. 
It is easy to see that after the softmax function application, neural network outputs become non-negative and their sum equals unity. This makes model's assumptions about the sample class mutually exclusive.
So, overall, we stated a function $f(x, \theta)$ where $\theta$ is a set of neural network parameters, which domain is rectangular representation matrices (of any complexity) and the set of values is 36-dimensional simplex.
\subsection{Training}
Essentially, training of the neural networks is adjusting their parameters (in our case $W_f$, $W_i$, $W_o$, $W_c$, $b_f$, $b_i$, $b_o$ and $b_c$ along with their counterparts for another direction and $W_{fc}$) based on empirical data. A common procedure for this adjustment which we used is called the \textbf{stochastic gradient descent}. In this subsection we're going to briefly describe this procedure. Namely, a set of training examples along with their known labels (this set is often called a \textbf{batch}, or \textbf{mini-batch}) 
$$\{(x_j, y_j), i \in (1...N)\}$$
is randomly sampled from the fixed training set, or, as in our case, obtained via stochastic generating procedure. In this context, $N$ is called a \textbf{batch size}. Then, for each training sample in the batch we evaluate the neural network output (see subsection \ref{architecture})
$$
\widetilde{y}_i = f(x_i, \theta).
$$

After that, we evaluate a \textbf{loss function}, which in our case is \textbf{categorical cross entropy}
$$
L = - \sum_{j=1}^{N}\sum_{i=1}^{36} I(y_j,i) log(\widetilde{y}_{ij})
$$
where $I(y_j,i)$ is an indicator function which is equals unity if  sample with index $j$ belongs to the class $i$ and zero otherwise and $\widetilde{y}_{ij}$ is corresponding coordinate of the neural network output for that sample. One can see that this function is lower when the neural network prediction are closer to the correct ones and equals zero when correct classes are predicted with highest possible confidence.
In order to adjust neural network parameters, a \textbf{gradient step}
$$
\theta \longrightarrow \theta - \lambda \nabla_{\theta}L,
$$
is performed, where $\lambda \in \mathbf{R}^{+}$ is a \textbf{learning rate}.
The gradient of the loss function with respect to the neural network is evaluated via backpropagation algorithm.

\section{Experiments}

\subsection{Data}
The nature of our task provides a vast amount of the training data. Though, obtaining sufficiently many diagrams with the known class is non-trivial. We prepare each batch of the training data for our model in the following way: uniformly sample with repetitions a set of 2048 (which is our batch size) 'classic' knot diagrams, perform randomly chosen stabilization moves until each diagram complexity reaches value uniformly sampled between 25 and 30 for each batch. Then 1000 commutation moves (including external)  are performed.  Validation set  consists of the 10000 diagrams of complexity 35 prepared with the same protocol and kept unchanged. 

\subsection{Importance of the external switch}
It turns out that set of moves which were picked to both train and test diagrams generation procedure has crucial impact on the neural network performance. We have conducted two experiments: in the first diagrams were perturbed using internal edge switching moves solely. In the second one, we have used the external switches as well. Here we will call these experiments internal and external case.

\begin{figure}[H]
  \center
  \includegraphics[width=0.99\linewidth]{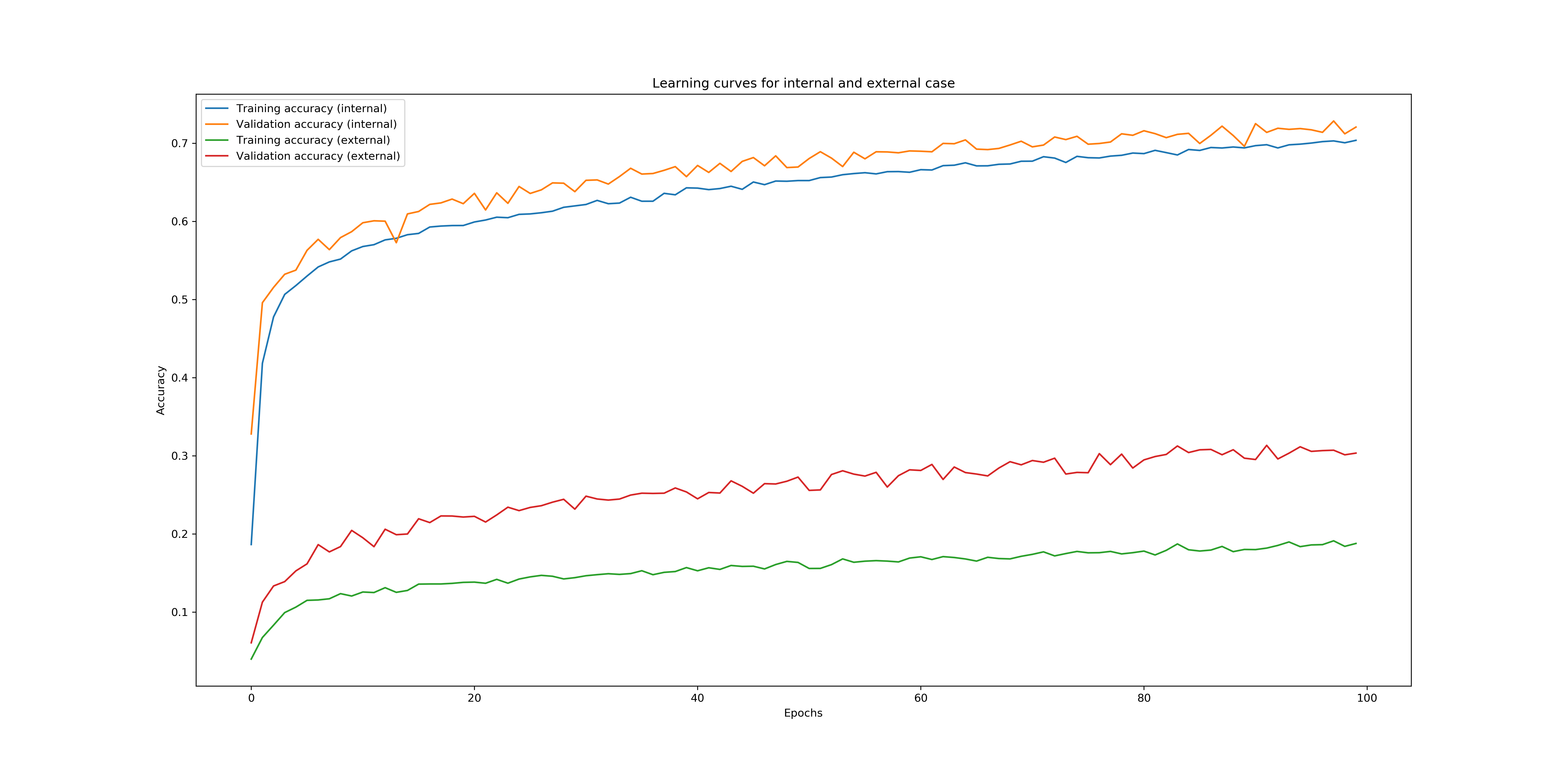}
    \caption{Learning curves for internal and external case}
  \label{fig:learning_curves_int_ext}
\end{figure} 

From the Fig.\ref{fig:learning_curves_int_ext} one can easily see that incorporation of external switch to the move set brings the complexity of the task to a new level. With exactly the same network architecture validation accuracy approaches much lower value. 

\subsection{Training and hyperparameters}
We trained our model for 300000 steps with Adam optimizer with learning rate 0.0001. Batch size was set to 2048. Since we do not have a fixed training set, we consider 100 gradient steps as an epoch. On the Fig.\ref{fig:learning_curves_int_ext} only first 100 epochs were considered. Full learning curves may be found in the section \ref{section: ext_results}.

\section{Results}
Our model managed to achieve significantly better accuracy than $100\%*(1/36) = 2.78 \%$ random guessing baseline. This accuracy is expected for a dummy algorithm which assigns knot classes randomly.
Speaking of the unknots in validation set, precision and recall  achieved 0.84 and 0.95 respectively. Full classification report, along with definitions of precision and recall may be found in the section \ref{section: ext_results}. One also may find interesting the confusion matrix for the validation set. The $(i,j)$ element of this matrix shows how many times the model considered the sample which belongs to the class $i$ as a sample belonging to the class $j$. Diagonal elements corrrespond to correct prediction, so in the ideal case this matrix is diagonal. In order to enhance its informativeness, values were binarized - elements greater than 40 showed with white pixels, everything else is black. One can see that besides mistakes we consider sporadic the model tends to confuse knot $5_1$ with $3_1$ and $7_1$.

\subsection{Dependence of accuracy on the number of complication steps}
We are getting samples by complicating the diagrams from finite set of 'classic' diagrams of knots with minimal number of intersections. If we don't complicate them enough, say, just learn on initial diagrams, they would be just memorized by large model. If we complicate them too much, the task becomes too hard for the model to perform. In order to illustrate dependence of prediction accuracy on number of complication steps we conduct the following experiment: we take a several test sets of different complexity, perform  500 complication steps  and then complicate our diagrams further measuring accuracy on the way  each 100 steps. 
\begin{figure}[H]
  \center
  \includegraphics[width=0.9\linewidth]{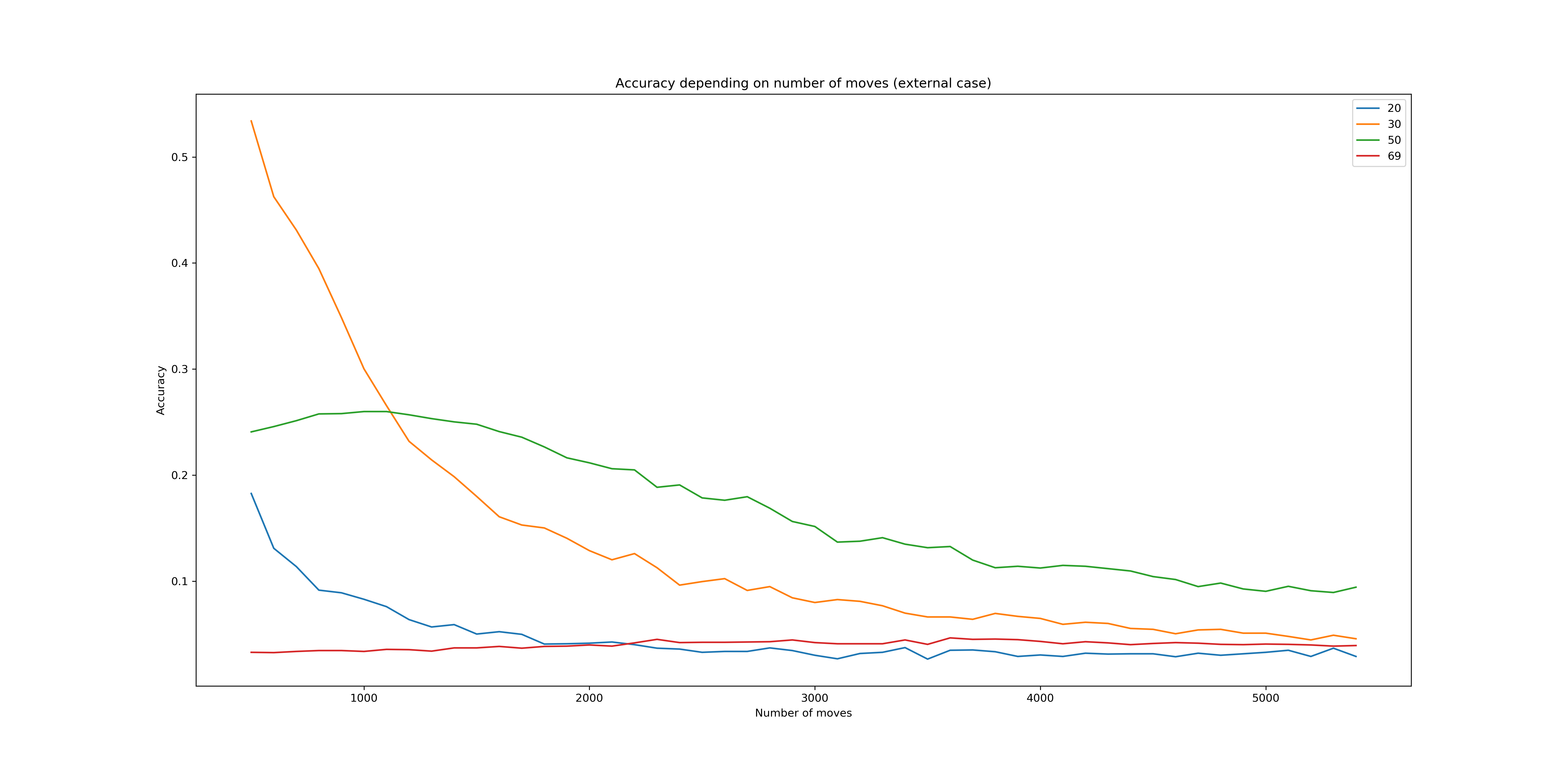}
  \caption{Accuracy with respect to the number of moves on the test set(external case)}
  \label{fig:acc_ent_external}
\end{figure} 
It turns out that levels of 15, 70 and 80 turned out to be the hardest cases for the models which makes sense since these kind of diagrams are the farthest from the training set. Complexity levels closer to the training ones (between 25 and 30) show the highest test set accuracy. Moreover, we can see that accuracy decreases for small diagrams and slightly increases for the large ones for some time then also decreases with their further complication. One may look at Fig. \ref{fig:heatmap_external} to explore dependence of accuracy regarding two factors - diagram complexity and number of entangling steps.

\subsection{Test-time augmentation}
Since  we have easy ways to perturb the diagram preserving its class, we find test-time augmentation a plausible way to increase the accuracy of predicting diagrams' classes.
To illustrate that, we conduct this kind of experiment: prepare test set of 3600 diagrams of complexity 35, and then slowly (1 move at a time) perturb it, compute the model output on changed diagrams (let's denote these predictions on the step $j$ as $y_j$) and then average the predictions. 
$$
\hat{y_j} = \dfrac{1}{j}\sum_{k=1}^{j}y_k
$$

The following plot shows that prediction accuracy benefits from adding more prediction steps even though accuracy of individual predictions may degrade.
\begin{figure}[H]
  \center
  \includegraphics[width=0.9\linewidth]{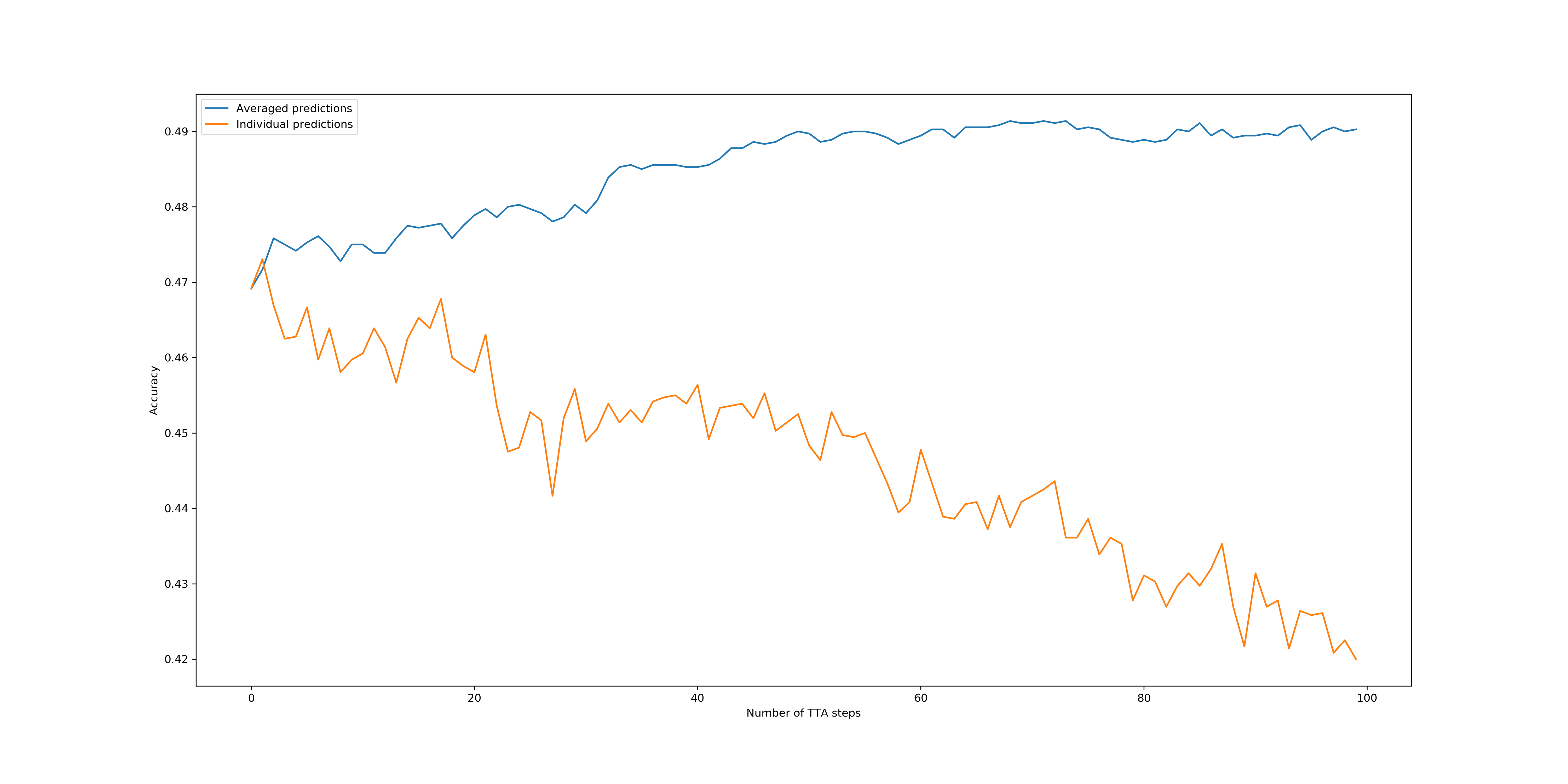}
  \caption{Result of the test-time augmentation (external case).}
  \label{fig:tta_external}
\end{figure} 
We believe that the same technique can be applied to the approach used in \cite{Vand} to enhance the prediction accuracy. There are some label-preserving perturbations of the knot in their parametrization such as rigid rotations of the whole polymer, which could be applied for the test-time augmentation.

\section{Discussion}
We've collected some evidence of efficacy of deep learning-based models to prediction of invariants of topological objects. There are many space for improvements left in terms of appropriate neural network architecture and hyperparameters choice, diagram parametrization and so on. In this paper we choose pretty simple neural network in order to make the whole process transparent to the readers with mathematical background. Our model didn't show any signs of overfitting thus might be trained longer with expectation of better performance. The knowledge of what makes it hard to correctly predict a knot diagram class is still to be found. There are infinitely many knot diagrams, as well as the classes, and designing the training and testing set, and the very definition of generalization in such case is a fundamental question. Definition of the arc diagram complexity as number of vertical edges does not always goes well with the experimental data: smaller diagrams sometimes are harder to classify by neural network.

Another deep learning approach seem promising: simplifying rectangular knot diagrams performing moves chosen by neural network.

\section{Appendix A}
\label{section: perm}
This appendix contains description of Dynnikov moves in terms of permutation representation. 
\subsection{Transition to the dual diagram}
Basically, any Dynnikov moves can be performed on both vertical and horizontal edges. So, we must find a way to switch between permutation representation in terms of vertical and horizontal edges. Let's consider a permutation representation in terms of horizontal edges:
$$
\begin{pmatrix}
3 & 8 & 9 & 1 & 6 & 2 & 4 & 5 & 7 \\
1 & 2 & 6 & 5 & 3 & 4 & 7 & 8 & 9
\end{pmatrix}.
$$
It contains horizontal coordinates of consecutive edges. Let's add corresponding vertical coordinates for both vertices in each edge as additional rows (the second and the fourth):
$$
\begin{pmatrix}
3 & 8 & 9 & 1 & 6 & 2 & 4 & 5 & 7 \\
1 & 2 & 3 & 4 & 5 & 6 & 7 & 8 & 9 \\
1 & 2 & 6 & 5 & 3 & 4 & 7 & 8 & 9 \\
1 & 2 & 3 & 4 & 5 & 6 & 7 & 8 & 9
\end{pmatrix}.
$$
The first two rows contain coordinates for the $X$ vertices, and the second - for the $O$ ones. Now let's sort $X$ and $O$ vertices (separately) according to their horizontal coordinates:
$$
\begin{pmatrix}
1 & 2 & 3 & 4 & 5 & 6 & 7 & 8 & 9 \\
4 & 6 & 1 & 7 & 8 & 5 & 9 & 2 & 3\\
1 & 2 & 3 & 4 & 5 & 6 & 7 & 8 & 9 \\
1 & 2 & 5 & 6 & 4 & 3 & 7 & 8 & 9
\end{pmatrix}.
$$
Now each pair of $X$ and $O$ vertices which belong to the same column of the obtained matrix has the same horizontal coordinates which, by definition of the rectangular diagram, makes them belonging to the same vertical edge. So, in order to obtain the permutation representation we just discard the rows with horizontal coordinates:
$$
\begin{pmatrix}
4 & 6 & 1 & 7 & 8 & 5 & 9 & 2 & 3\\
1 & 2 & 5 & 6 & 4 & 3 & 7 & 8 & 9
\end{pmatrix}.
$$
\subsection{Switch moves}
Switch moves are easily interpreted in terms of the permutation notation. Let's recall our diagram from the Fig.[\ref{fig:init_diag}] and its representation
$$
\begin{pmatrix}
3 & 8 & 9 & 1 & 6 & 2 & 4 & 5 & 7 \\
1 & 2 & 6 & 5 & 3 & 4 & 7 & 8 & 9
\end{pmatrix}.
$$
Let's consider the diagram which is one move away from the initial one:

\begin{figure}[H]
  \center
  \includegraphics[width=0.9\linewidth]{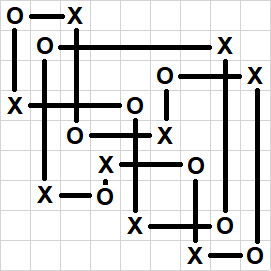}
  \caption{An example of the diagram one internal switch move avay from \ref{fig:init_diag}}
  \label{fig:inter_switch_diag}.
\end{figure}
Its permutation representation is
$$
\begin{pmatrix}
3 & 8 & 9 & 1 & 6 & 4 & 2 & 5 & 7 \\
1 & 2 & 6 & 5 & 3 & 7 & 4 & 8 & 9
\end{pmatrix}.
$$
Since switching of horizontal edges affects only their order, not coordinates of their ends, it's only a matter of order of columns in the corresponding permutation matrix. 

\section{Appendix B}
\label{section: diag_examples}
This appendix contains rectangular diagrams used for the neural network training.
\begin{figure}[H]
  \center
  \includegraphics[width=0.9\linewidth]{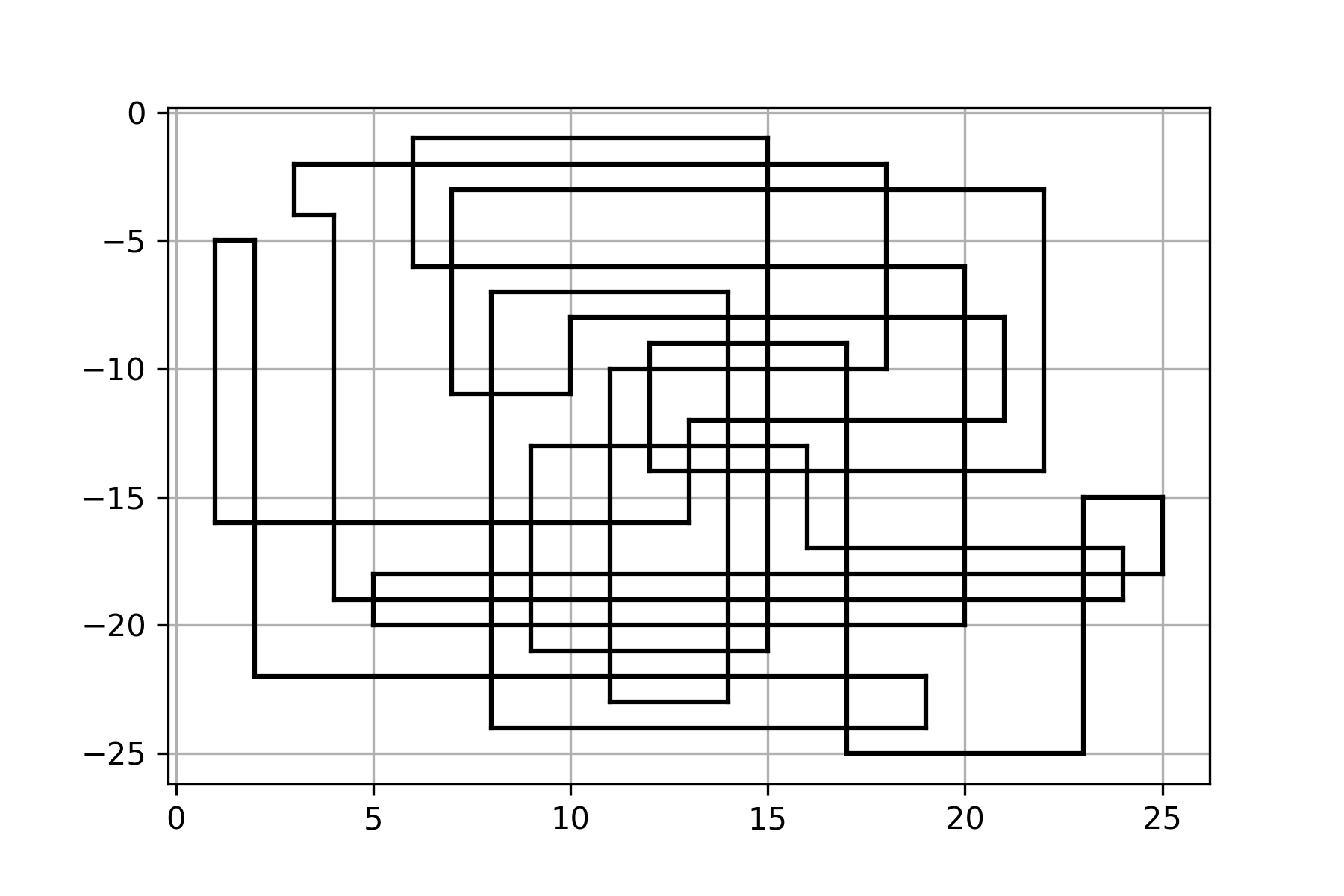}
  \caption{Unknot diagram of complexity 25, 1000 moves}
\end{figure}

\begin{figure}[H]
  \center
  \includegraphics[width=0.9\linewidth]{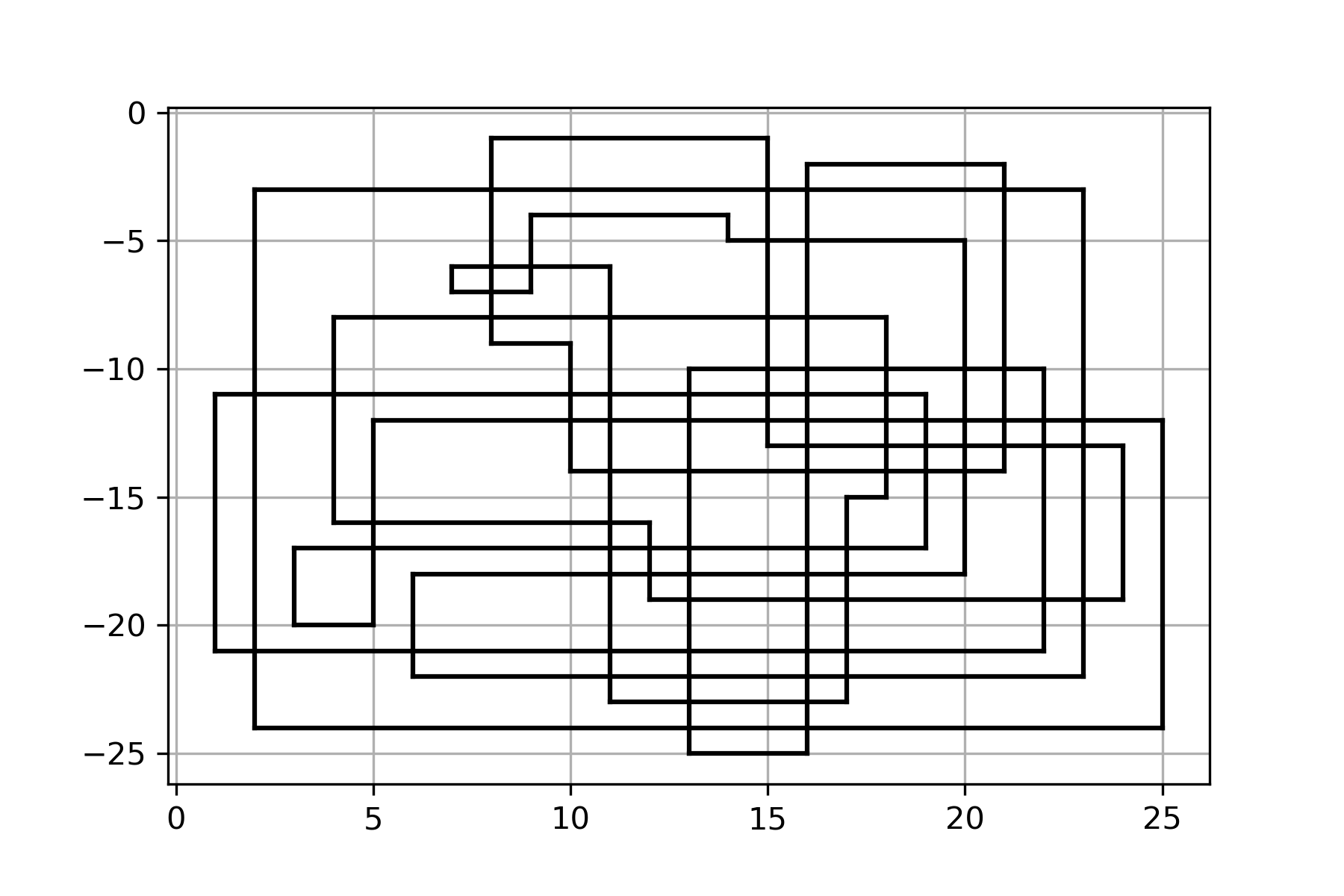}
  \caption{Knot diagram of complexity 25, 1000 moves}
\end{figure} 

\section{Appendix C} 
This appendix contauins illustrations regarding results of the neural network training and evaluation.
\label{section: ext_results}

\begin{figure}[H]
  \label{fig:learning_curve_acc}
  \center
  \includegraphics[width=0.9\linewidth]{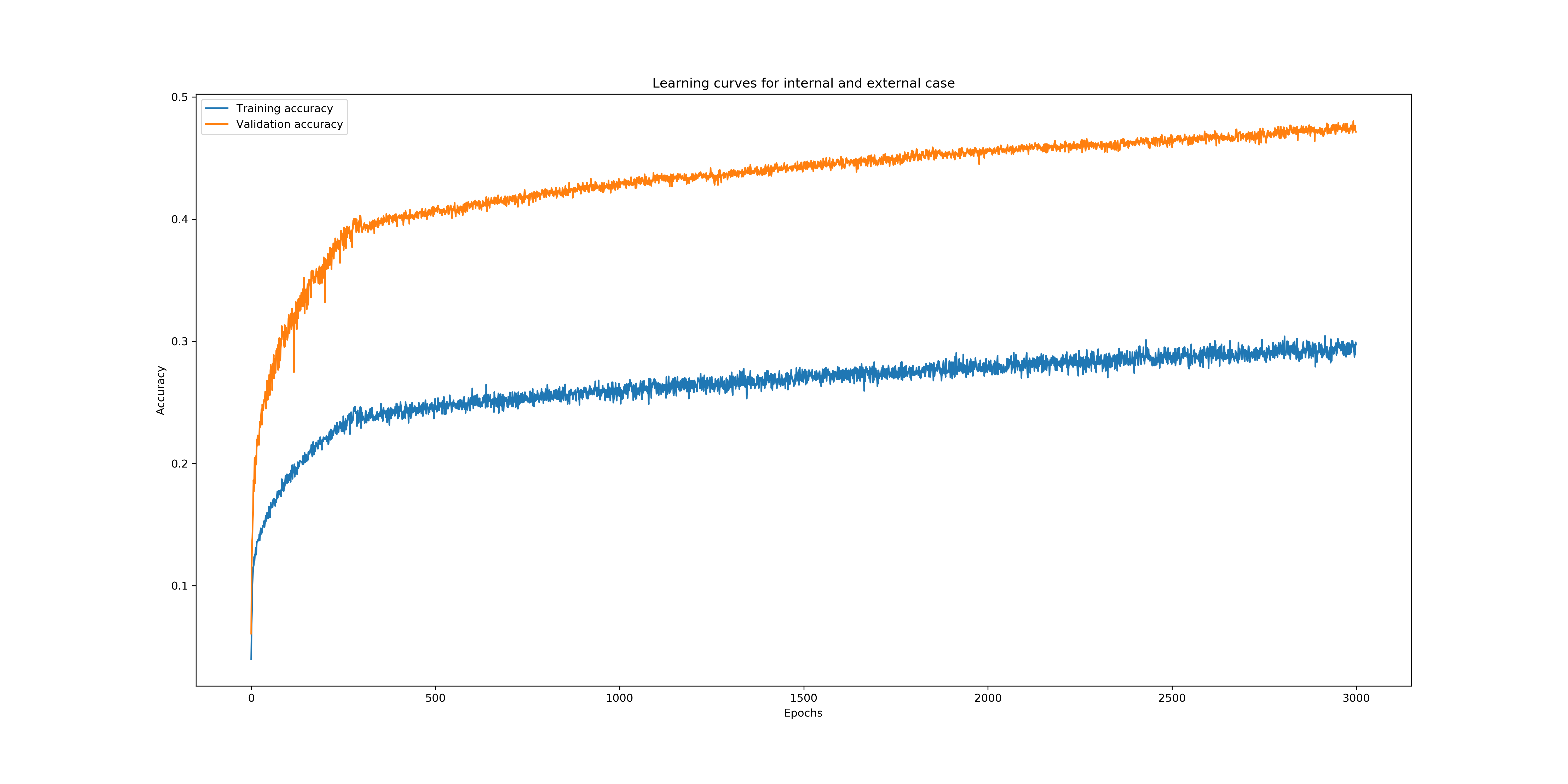}
  \caption{Learning curve (accuracy) }
\end{figure}

\begin{figure}[H]
  \label{fig:learning_curve_loss}
  \center
  \includegraphics[width=0.9\linewidth]{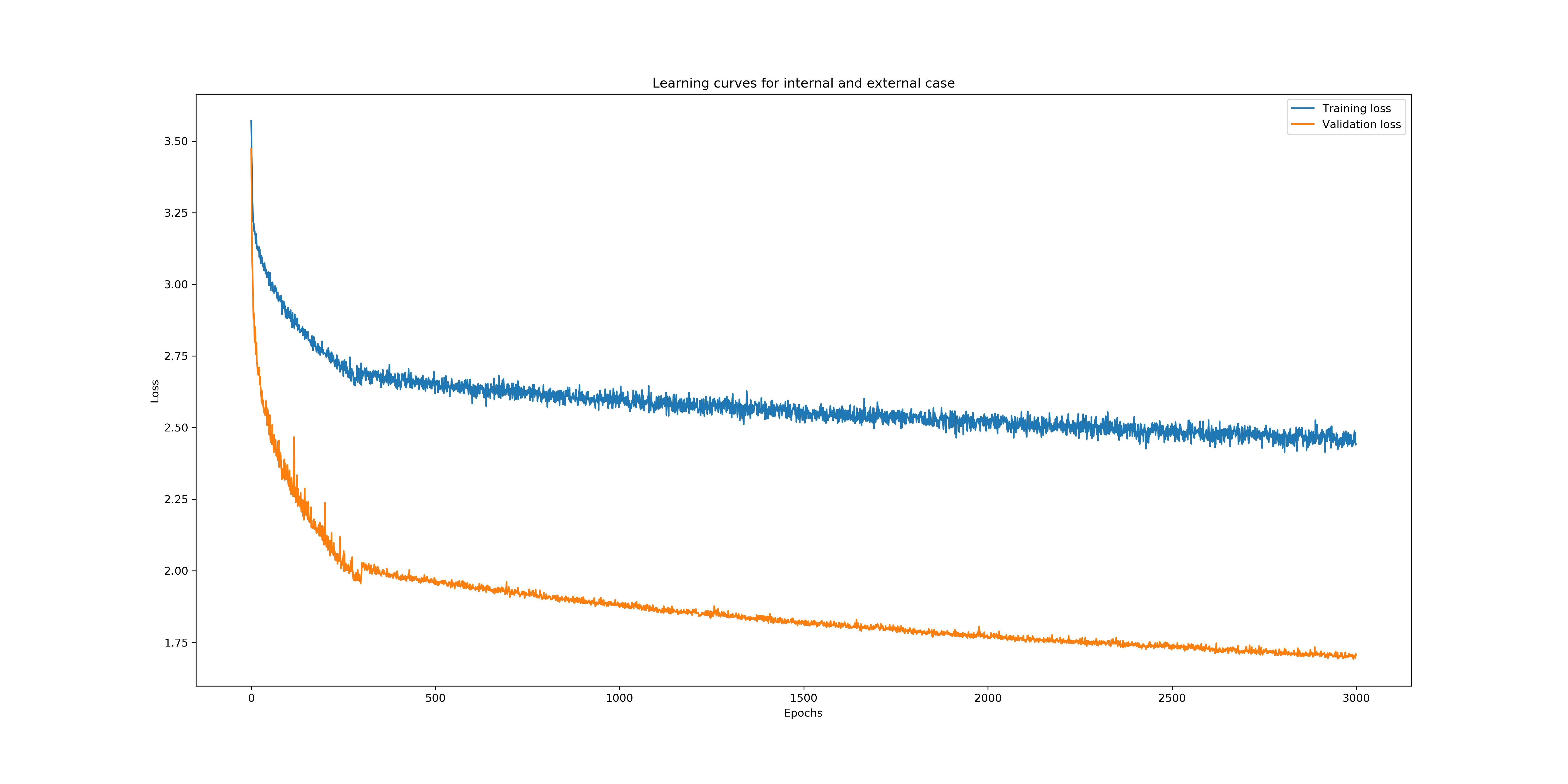}
  \caption{Learning curve (loss) }
\end{figure} 

\begin{figure}[H]
  \center
  \includegraphics[width=0.9\linewidth]{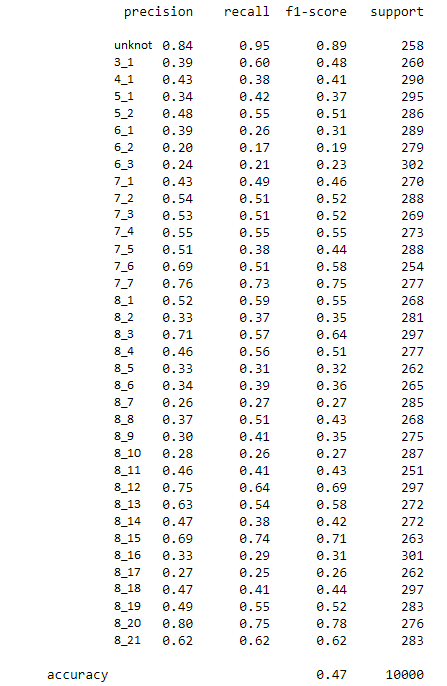}
  \caption{Classification report for all considered knot classes}
\end{figure}
Let us explain the meaning of each column in the classification report. For each class, precision and recall are evaluated with the following formulas:
$$
precision = \dfrac{TP}{TP+FP},
$$
$$
recall = \dfrac{TP}{TP+FN},
$$
where $TP$ is number of samples in the evaluation set correctly assigned to the class (the model assigned the sample its true class), $FP$ is number of the samples of another classes which were assigned by the model to considered class, and $FN$ is number of the samples with considered class assigned to other classes. In other words, for each class, precision is a fraction of correct predictions among samples were assigned to this the class by the model, and recall is a fraction of correct predictions among samples which truly belong to this class.
Then, f1-score is a harmonic mean of precision and recall
$$
f1-score = \dfrac{2*precision*recall}{precision+recall}
$$
and is aimed to have a single number for each to measure how good this class is handled by the model. A support is just a number of samples in the validation set which belong to each class. An accuracy is a fraction of correctly classified samples.

\begin{figure}[H]
  \center
  \includegraphics[width=0.9\linewidth]{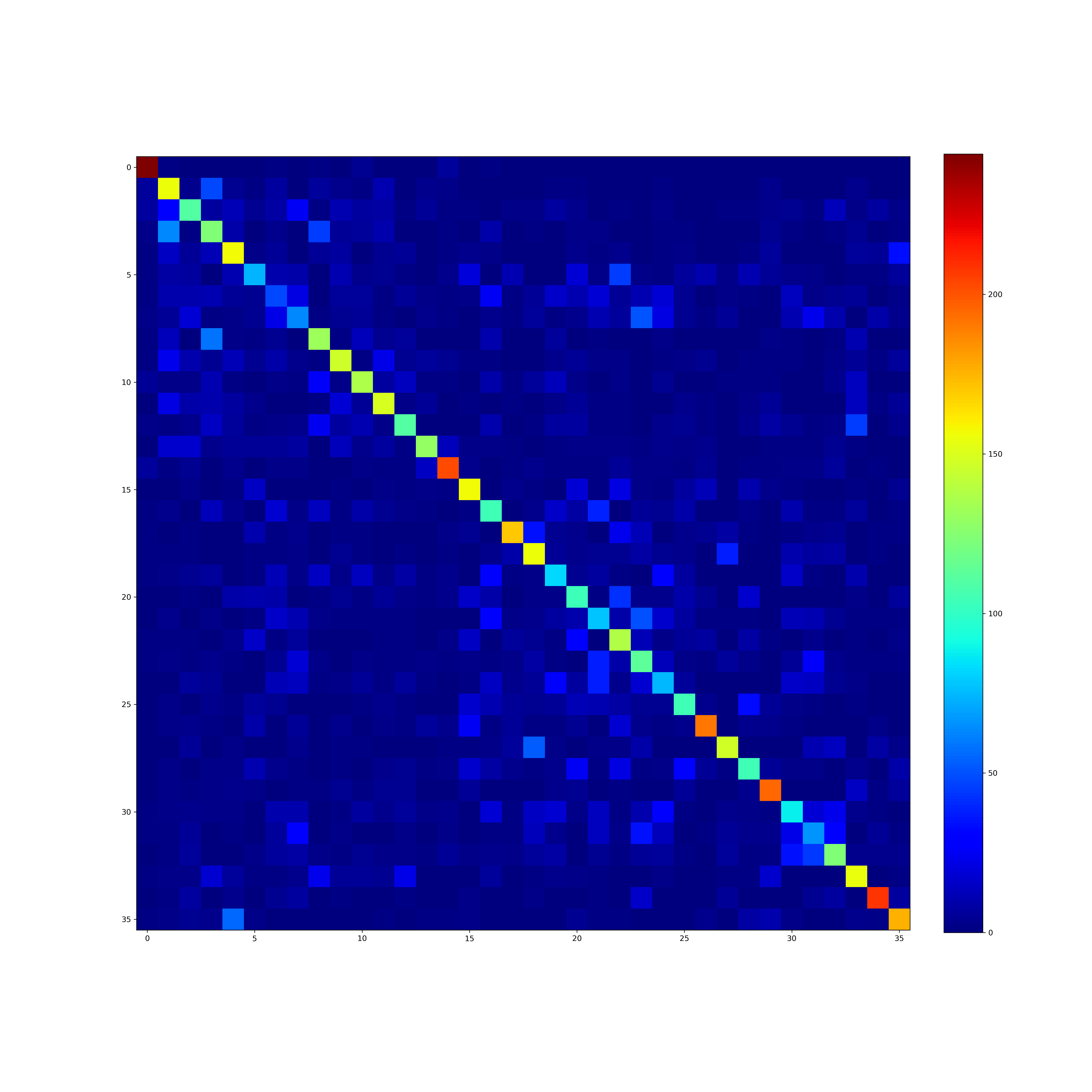}
  \caption{Confusion matrix for the test set}
  
\end{figure}

\begin{figure}[H]
  \center
  \includegraphics[width=0.9\linewidth]{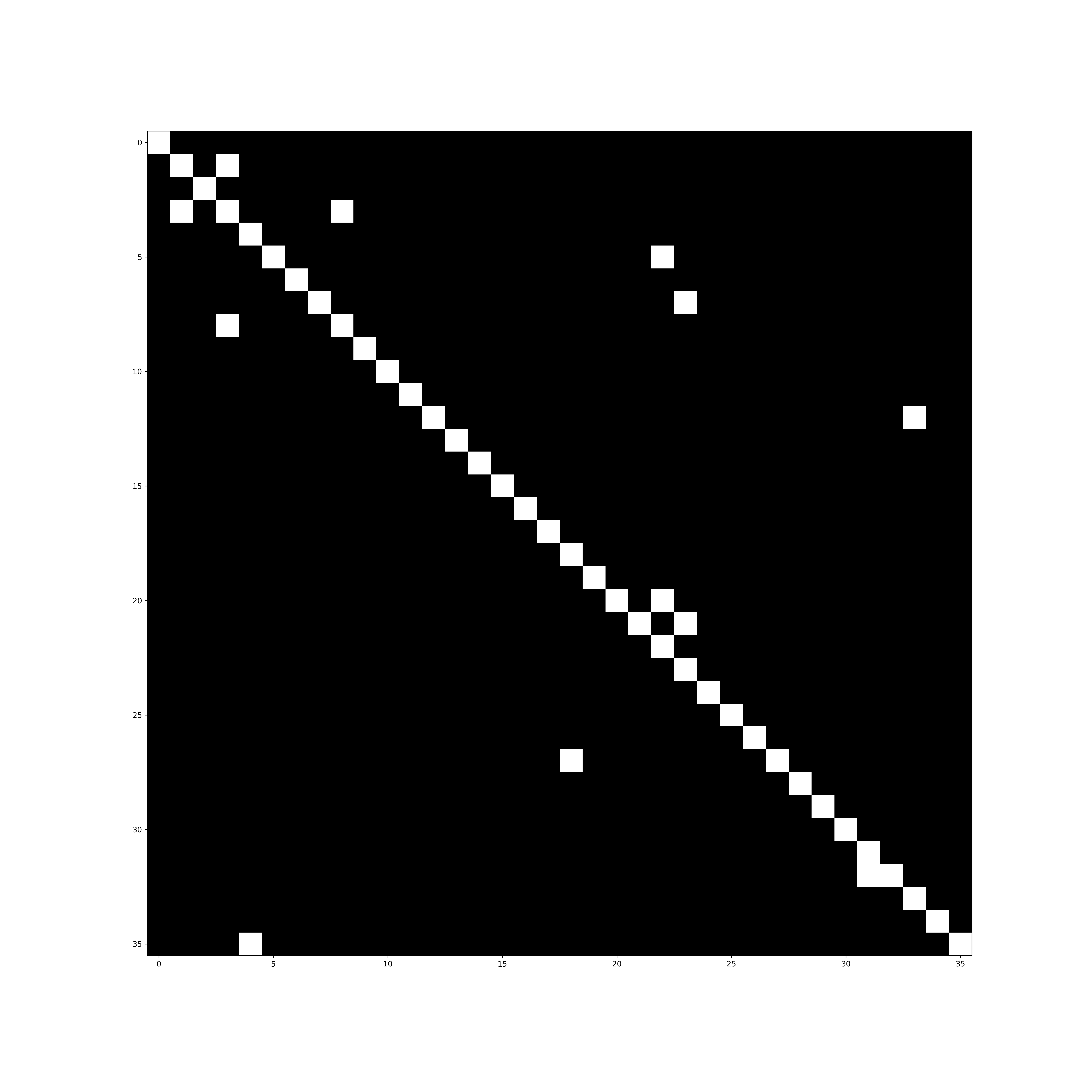}
  \caption{Binarized confusion matrix(the threshold is 40)}
  
\end{figure}

\begin{figure}[H]
  \center
  \includegraphics[width=0.9\linewidth]{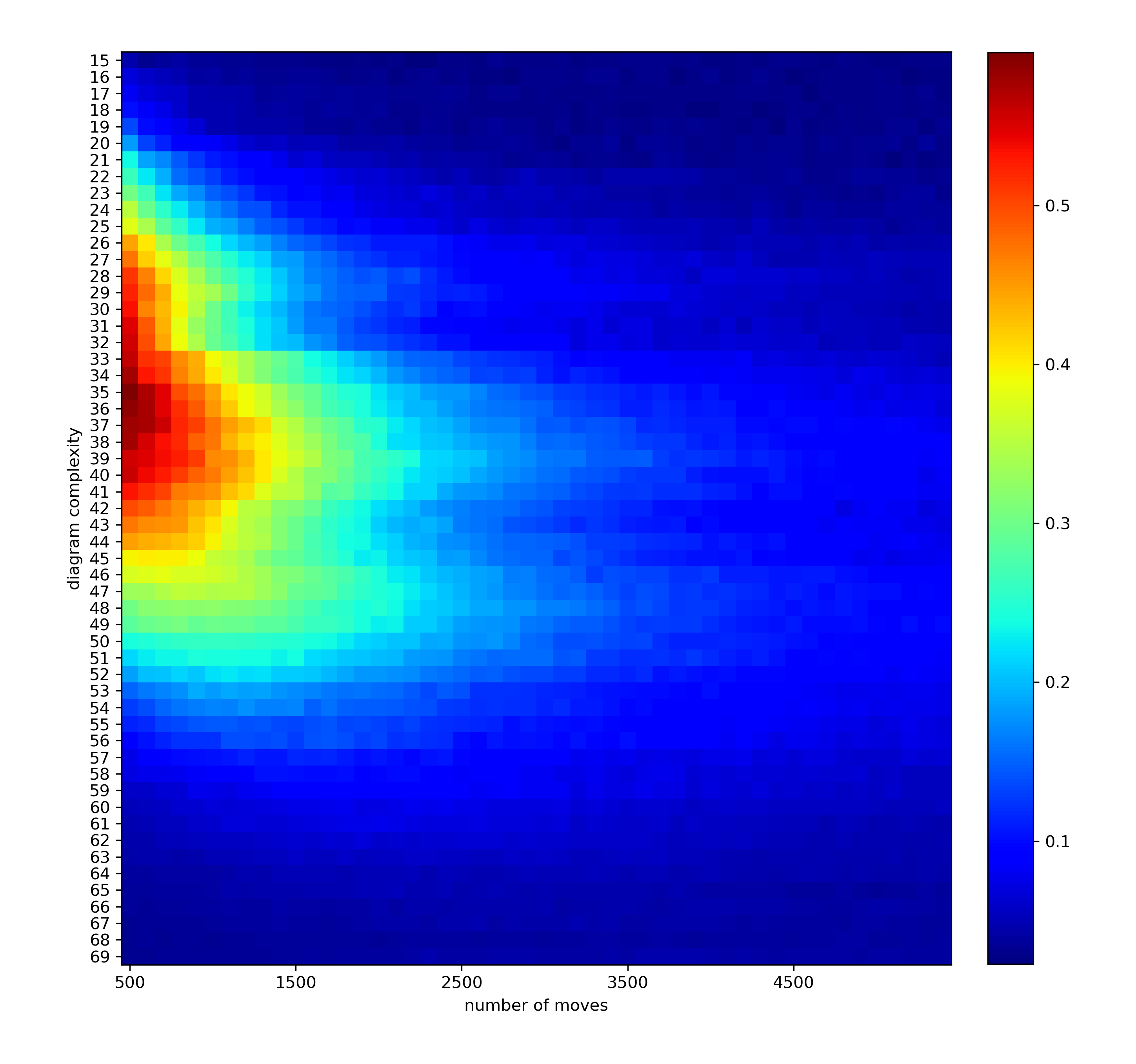}
  \caption{Accuracy heatmap(external case)}
  \label{fig:heatmap_external}
\end{figure}

\bibliographystyle{unsrt}  


\end{document}